\documentclass{elsart3-1}

\usepackage{latexsym}
\usepackage{mathtools}
\usepackage{amssymb}
\usepackage{amsmath}
\usepackage[dvips]{graphicx}
\newlength{\cqfd}
\usepackage[english,francais]{babel}
\newcommand{\into}{\int_{\Omega}}
\newcommand{\he}{h_\varepsilon}
\newcommand{\ue}{\bu_\varepsilon}
\newcommand{\ep}{\varepsilon}
\newcommand{\ve}{v_\varepsilon}

\newcommand{\bu}{\bar{u}}
\newtheorem{theorem}{Theorem}[section]

\newtheorem{e-proposition}[theorem]{Proposition}

\newtheorem{e-definition}[theorem]{Definition\rm}

\renewcommand{\theequation}{\arabic{equation}}
\def\og{\leavevmode\raise.3ex\hbox{$\scriptscriptstyle\langle\!\langle$~}}
\def\fg{\leavevmode\raise.3ex\hbox{~$\!\scriptscriptstyle\,\rangle\!\rangle$}}

\journal{the Acad\'emie des sciences}
\begin{document}
\centerline{}
\begin{frontmatter}


\selectlanguage{english}
\title{BD Entropy and   \\
        Bernis-Friedman Entropy}


\selectlanguage{english}
\author[authorlabel1]{Didier Bresch\thanksref{lab1}},
\ead{Didier.Bresch@univ-smb.fr}
\author[authorlabel2]{Mathieu Colin\thanksref{lab2}},
\ead{Mathieu.Colin@math.u-bordeaux.fr}
\author[authorlabel1]{Khawla Msheik\thanksref{lab1}},
\ead{khawla.msheik@univ-smb.fr}
\author[authorlabel3]{Pascal Noble\thanksref{lab3}},
\ead{pascal.noble@math.univ-toulouse.fr}
\author{Xi Song}
\ead{xi.song@hotmail.com}
\thanks[lab1]{Research of D.B. and K.M. was partially supported by the ANR project ANR-16-CE06-0011-02 FRAISE and  ANR-17-CE08-0030 ViSCAP. D.B. want to thank G. Kitavtsev for his comments and suggestions.}
\thanks[lab2]{Research of M.C. was partially supported by the ANR project ANR-17-CE08-0030 ViSCAP}
\thanks[lab3]{Research of P.N. was partially supported by the ANR project BoND ANR-13-BS01-0009-01 .}
\address[authorlabel1]{LAMA UMR5127 CNRS, Universit\'e Grenoble Alpes, Universit\'e Savoie Mont-Blanc, 73376 Le Bourget du Lac, France}
\address[authorlabel2]{Equipe INRIA CARDAMOM , IMB Equipes EDP, 351 cours de la  lib\'eration, 33405 Talence, France}
\address[authorlabel3]{IMT, INSA Toulouse, 135 avenue de Rangueil, 31077 Toulouse Cedex 9, France}


\medskip
\begin{center}
{\small Received *****; accepted after revision +++++\\
Presented by }
\end{center}

\begin{abstract}
\selectlanguage{english}
In this note, we  propose in the full generality a  link between the BD entropy introduced by 
D.~Bresch and B.~Desjardins for the viscous shallow-water equations and the Bernis-Friedman 
(called BF) dissipative entropy  introduced to study the lubrications equations.
Different dissipative entropies are obtained playing with the drag terms on the viscous shallow
water equations. It helps for instance to prove global existence of 
nonnegative weak solutions for the lubrication equations starting from the global existence
of nonnegative weak solutions for appropriate viscous shallow-water equations.

\medskip

\noindent{\bf R\'esum\'e} \vskip 0.5\baselineskip \noindent
{\bf BD entropie et entropie de Bernis-Friedman. }
 Dans cette note, on propose un lien g\'en\'eral entre la BD entropie introduite par
 D. Bresch et B. Desjardins pour les \'equations de Saint-Venant visqueuses et l'entropie 
 dissipative de Bernis-Friedman  (not\'ee BF) introduite pour \'etudier les \'equations de lubrifications. Diff\'erentes entropies dissipatives sont obtenues suivant le choix des termes de train\'ee sur Saint-Venant visqueux.
   Ce lien entre ces deux outils math\'ematiques  aide par exemple \`a prouver l'existence de solutions faibles  positives pour les \'equations de lubrification en partant de l'existence de solutions faibles 
 positives pour des \'equations  de Saint-Venant visqueuses bien choisies.

\end{abstract}
\end{frontmatter}

\section*{Version fran\c{c}aise abr\'eg\'ee}
Dans cette note,   on propose un lien g\'en\'eral entre la BD entropie introduite par
 D.~Bresch et B.~Desjardins pour les \'equations de Saint-Venant visqueuses (voir \cite{BD} et \cite{BrDe2}) et l'entropie  dissipative de Bernis-Friedman (que l'on notera BF) introduite (voir \cite{BF})  pour \'etudier les \'equations de lubrifications. Diff\'erentes entropies dissipatives sont obtenues suivant le choix des termes de train\'ee sur Saint-Venant visqueux g\'en\'eralisant ainsi quelques travaux importants comme \cite{KiLaNi}, \cite{FoKiTa}.  Ce lien entre les deux outils importants que sont la BD entropie et la BF entropie permet par exemple de construire des solutions faibles du mod\`ele
 de lubrification \`a partir de solutions faibles du mod\`ele de Saint-Venant.   Il permet
 \'egalement d'obtenir certains r\'esultats sur les \'equations de Saint-Venant
 en s'inspirant des r\'esultats \'etablis sur les \'equations de lubrification qui ont  \'et\'e beaucoup
 plus \'etudi\'ees historiquement.  Le syst\`eme de lubrification s'\'ecrit par exemple
 $$\partial_t h + \partial_x( \frac{1}{\alpha W_e} h^n \partial_x^3 h- \frac{1}{\alpha Fr^2} h^{m-1} \partial_x h) = 0
 $$
 et le mod\`ele de Saint--Venant associ\'e
\begin{equation}\label{swgg}
\begin{cases}
\partial_t h_\varepsilon+\partial_x (h_\varepsilon\bu_\varepsilon)=0,\\
\varepsilon \Bigg(\partial_t (h_\varepsilon \bu_\varepsilon)+\partial_x(h_\varepsilon \bu_\varepsilon^2) \Bigg)  +\frac{1}{Fr^2}h_\varepsilon^{\beta} \partial_x(h_\varepsilon)=\varepsilon \Bigg(\frac{4}{{ R_e}}\partial_x(h_\varepsilon\partial_x \bu_\varepsilon) \Bigg) + \frac{1}{{ W_e}}h_\varepsilon\partial_x^3 h_\varepsilon 
- \alpha \frac{h_\varepsilon^2 \overline u_\varepsilon }{\he^n},
\end{cases}
\end{equation}
o\`u $\beta+n = m$: voir par exemple \cite{Bert1}, \cite{Bert2} ou \cite{BD}. Nous discuterons des deux outils importants que sont la BF entropie et l'entropie dissipative  d\^ue \`a Bernis-Friedman (BF). Nous expliquerons l'int\'er\^et qu'il y a de  mettre en exergue une telle relation. On peut par exemple \'etudier ces syst\`emes et montrer l'existence globale de solutions faibles du mod\`ele de lubrification en partant de solutions faibles positives du mod\`ele de Saint-Venant associ\'e.
   On peut \'egalement  consid\'erer des syst\`emes avec des termes non-locaux en s'inspirant de r\'esultats r\'ecents, voir par exemple \cite{ImMe}, \cite{CaWrZa}.

\selectlanguage{english}

\section{Introduction: Lubrication systems and viscous shallow-water equations with drag terms.}
\label{sec2}
\noindent
      In this section, we present the formal link between two key tools respectively for lubrication
system by Bernis-Friedman and for shallow-water equations by Bresch-Desjardins. We first
start by presenting the two quantities and their link on a simple example and then we explain
how to get relations in a more general case. Our calculations remain at this stage only formal.
We assume solutions are regular enough.

\subsection{A lubrication system: energy estimate and Bernis-Friedman (BF) dissipative entropy}
 In a one dimensional periodic domain $\Omega $,  consider the following thin-film equation, also known as lubrication equation
\begin{equation}\label{lub1}
\partial_t h+\partial_x ( \frac{1}{\alpha W_e}F(h)\partial_x^3 h-\frac{1}{\alpha Fr^2}F(h)\partial_x h)=0.
\end{equation}
 We couple this equation with the initial condition
$$ h(x,0)= h_0(x) \quad \mbox{in } \Omega$$
System (\ref{lub1}) can be rewritten equivalently as a gradient flow system
\begin{equation}\label{lub2}
\begin{cases}
\partial_t h+\partial_x (h u)=0,\\
h u = \frac{1}{\alpha W_e}F(h)\partial_x^3 h-\frac{1}{\alpha Fr^2}F(h)\partial_x h.
\end{cases}
\end{equation}
The corresponding energy is given, for all $t\in (0,T)$,  by
\begin{equation} \label{en1}
\int_{0}^{t} \into \frac{\alpha h^2 u^2}{F(h)}~dx~dt + \frac{1}{2}\into   \frac {h(x,t) ^2}{Fr^2} + \frac{(\partial_x h(x,t)) ^2}{W_e} ~dx =  \frac{1}{2}\into \frac{h_0(x) ^2}{Fr^2} + \frac{(\partial_x h_0(x)) ^2}{W_e} ~dx.
\end{equation}
In their paper \cite{BF}, Bernis and Friedman proved the existence of a weak solution for a higher order nonlinear degenerate parabolic equations and suggested a new entropy inequality- referred to by BF entropy- which provides additional estimates serving for increasing the regularity of the weak solution obtained. As for our problem, we adapt the same methodology to obtain the BF entropy of the general lubrication model stated above. Indeed, define the functionals
$$g_{\ep}(s)= -\int_{s}^{A} \frac{1}{F(r)+ \ep}\,dr, \qquad G_{\ep}(s)= -\int_{s}^{A} g_{\ep}(r) \,dr,$$
with $A$ being an integer such that $A\ge \max|h(x,t)|$.
According to Bernis and Friedman, we multiply (\ref{lub1}) by $G_0'(h)$, where $G_0= \lim\limits_{\ep \rightarrow 0} G_{\ep}$, we get the  BF dissipative entropy equality
\begin{equation}\label{bf}
\into G_0(h(x,T))~dx  + \int_{0}^{T} \into \frac{(\partial_x^2h)^2}{\alpha W_e} + \frac{(\partial_xh)^2}{\alpha Fr^2}~dx~dt =\into G_0(h_0(x))~dx  .
\end{equation}
 \subsection{A viscous shallow--water system: energy estimate and BD entropy}
 In a periodic domain $\Omega $, let us consider the viscous Shallow Water system with surface tensions and a drag term: 
\begin{equation}\label{swg}
\begin{cases}
\partial_t h_\varepsilon+\partial_x (h_\varepsilon\bu_\varepsilon)=0,\\
\varepsilon \Bigg(\partial_t (h_\varepsilon \bu_\varepsilon)+\partial_x(h_\varepsilon \bu_\varepsilon^2) \Bigg)  +\frac{h_\varepsilon \partial_x(h_\varepsilon)}{Fr^2}=\varepsilon \Bigg(\frac{4}{{ R_e}}\partial_x(h_\varepsilon\partial_x \bu_\varepsilon) \Bigg) + \frac{1}{{ W_e}}h_\varepsilon\partial_x^3 h_\varepsilon 
- \alpha \frac{h_\varepsilon^2 \overline u_\varepsilon }{F(h_\varepsilon)}.
\end{cases}
\end{equation}
The initial conditions are given by
$$ h_\varepsilon|_{t=0} = h_0^\varepsilon,
\qquad
(h_\varepsilon  \bu_\varepsilon) |_{t=0} = m_0^\varepsilon.
$$
$\alpha$ is a positive constant, $ {R_e}$, ${ W_e}$ and ${ Fr}$ are respectively the adimentional Reynold, Weber and Froude numbers. Note that the terms in the right-hand side of the momentum equation represent respectively the viscous term, the capillary term and the drag term. The energy equation corresponding to (\ref{swg}) is given, for all $t\in (0,T)$ by
\begin{align}\label{ensw}
& \left( \into\varepsilon \frac{\he(x,t) \bu_\varepsilon^2(x,t)}{2}+\frac{h_\varepsilon(x,t)^2}{2Fr^2}+\frac{(\partial_x h_\varepsilon(x,t)^2}{2W_e}\right) ~dx +\int_{0}^{t}\into\frac{4\varepsilon }{R_e}h_\varepsilon (\partial_x \bu_\varepsilon)^2+  \alpha \frac{h_\varepsilon^2 \overline u_\varepsilon^2}{F(h_\varepsilon)} ~dx~dt\\
 & \nonumber = \frac{1}{2}\into \ep \frac{(m_0^{\ep})^2}{h_0^{\ep}} + \frac{(h_0^{\ep})^2}{Fr^2} + \frac{(\partial_x h_0^{\ep})^2}{W_e}~dx.
\end{align}
As introduced in \cite{BD}, the BD entropy equality is obtained by deriving the mass equation in space and multiplying by $4\varepsilon/ R_e$ then summing with the momentum equation, and multiplying the sum by the artificial velocity $$ v_\varepsilon= \ue +  \frac{4}{R_e} \partial_x(log(\he)).$$
The BD entropy for the general system is given, for all $t\in (0,T)$, by
\begin{equation}
\begin{aligned}
&\frac{\ep}{2} \into \he(x,t) \ve(x,t)^2 ~dx  
+  \frac{1}{2}\bigg[ \into \frac{\he(x,t) ^2}{Fr^2} + \frac{(\partial_x \he(x,t)) ^2}{W_e} ~dx\Bigg] \\
&+ \alpha \int_{0}^{t}\into  \frac{\he^2 \ue ^2}{F(\he)} ~dx ~dt
+\frac{4 }{ R_e} \Bigg[ \into \int_{0}^{t} \frac{(\partial_x \he)^2}{ Fr^2}  
 +  \frac{(\partial_{xx} \he)^2}{ W_e} ~dx~dt + \alpha  \int_{0}^{t} \underbrace{\into  \frac{\he \ue }{F(\he)}\partial_x \he ~dx~dt}_{X}\Bigg]\\
&= \frac{\ep}{2} \into  \into h_0^{\ep}(x) \ve(x,0)^2 ~dx + \frac{1}{2} \into  \frac{ h_0^{\ep}(x)^2}{Fr^2} +\frac{(\partial_x h_0^{\ep}(x))^2}{W_e} ~dx .
\end{aligned}
\end{equation}
As for the term X, it can be rewritten as
\begin{align*}
X= \into \he \ue\frac{\partial_x \he }{F(\he)} \,dx &= \into \he \ue \big( \frac{d}{dx} \int_{A}^{\he} \frac{1}{F(y)}\,dy \big)\,dx \\
& = \into -\partial_x (\he \ue)  \big(  \int_{A}^{\he} \frac{1}{F(y)}\,dy \big)\,dx\\
&= \into \partial_t \he \big(  \int_{A}^{\he} \frac{1}{F(y)}\,dy \big) \,dx.
\end{align*}

\subsection{The link between the BD entropy and the BF dissipative entropy}
In view of the term X and $G_0$, we noticed the following
$$ G_0'(\he)= g_0(\he)= -\int_{\he}^{A} \frac{1}{F(r)}\,dr = \int_{A}^{\he} \frac{1s}{F(r)} ~dr . $$
Hence
\begin{align*} 
X = \into \partial_t\he G_0'(\he)\,dx = \frac{d}{dt} \into G_0(\he)\,dx.
\end{align*}
Thus, coupled with (\ref{ensw}), the BD entropy reads:
\begin{equation}
\begin{aligned}
&\frac{\ep}{2} \into \he(x,t) \ve(x,t)^2 ~dx - \frac{\ep}{2} \into \he(x,t) \ue(x,t)^2 ~dx  
-\frac{4\varepsilon }{R_e} \int_{0}^{t}\into h_\varepsilon (\partial_x \bu_\varepsilon)^2 \\
&
+\frac{4}{ R_e} \Bigg[ \into \int_{0}^{t} \frac{(\partial_x \he)^2}{ Fr^2}   
+  \frac{(\partial_{xx} \he)^2}{ W_e} ~dx~dt +  \alpha  \into  G_0(\he(x,t)) ~dx\Bigg]\\
&=-  \frac{1}{2}\into \ep \frac{(m_0^{\ep})^2}{h_0^{\ep}}+  \frac{\ep}{2}  \into h_0^{\ep}(x) \ve(x,0)^2 ~dx +   \frac{4\alpha }{R_e}\into  G_0(h_0^{\ep }) ~dx .
\end{aligned}
\end{equation}
  If we assume now that the couple $(\he,\ue)$, solution of (\ref{swg}), converges in a proper sense to $(h,u)$, then we find that the above BD entropy equality degenerates to the following inequality, which coincides with the dissipative BF-entropy of (\ref{lub1})
\begin{equation}
\begin{aligned}
 \into \int_{0}^{t} \frac{(\partial_x h)^2}{\alpha Fr^2}   +  \frac{(\partial_{xx} h)^2}{\alpha W_e} ~dx~dt +   \into  G_0(h(x,t)) ~dx =\into  G_0(h_0(x)) ~dx.
\end{aligned}
\end{equation}
Of course these computations are formal and have been written with equalities but they help
to understand that the BF entropy may be obtained from the BD entropy. This provides a way to
construct nonnegative solutions of lubrication equation from nonnegative solutions of the 
shallow water equation with appropriate drag terms. 
  Let us present below a general  computation with different surface tension and pressure term
showing the relation between the BD entropy and the BF dissipative entropy. 
 
\section{A general link between the BD entropy and the BF dissipative entropy}
\label{sec1}
In this part, we will consider the following fourth order lubrication approximation that has been studied in several papers, see for instance \cite{Bert1}, \cite{Bert2}:
\begin{equation}\label{lubg0}
\partial_t h+\partial_x ( \frac{1}{\alpha W_e}F(h)\partial_x^3 h-\frac{1}{\alpha Fr^2}D(h)\partial_x h)=0.
\end{equation}
In \cite{Bert1} for instance, the authors considered the above lubrication model with the following choice of $F$ and $D$: $F(h)= h^n$ and $D(h)=  h^{m-1}.$ Indeed, the lubrication equation becomes:
\begin{equation}\label{lubg}
\partial_t h+\partial_x ( \frac{1}{\alpha W_e}h^n\partial_x^3 h-\frac{1}{\alpha Fr^2}h^{m-1}\partial_x h)=0.
\end{equation}
 In fact, they proved the existence of a global in time nonnegative weak solution starting from nonnegative datum for all $n >0$, and $1<m<2$. In particular, The most critical case is the most significantly physical one when $n=3$ (moving contact line in a thin film). In this case, a distributional solution is proven to exist, where it becomes a strong positive solution in the infinite time limit. In this sequel, we will consider the choices of F and D stated above. Then, the BF entropy corresponding to the latter system is given by
\begin{equation}\label{bfg}
\into G_0(h(x,T)) ~dx + \int_{0}^{T} \into \frac{(\partial_x^2h)^2}{ We} +\frac{1}{ Fr^2} h^{m-n-1}(\partial_xh)^2 ~dx~dt= \into G_0(h_0(x)) ~dx.
\end{equation}
Consider herein the following Shallow water system with drag term corresponding to a weight $F(h)=h^n$
\begin{equation}\label{swgg}
\begin{cases}
\partial_t h_\varepsilon+\partial_x (h_\varepsilon\bu_\varepsilon)=0,\\
\varepsilon \Bigg(\partial_t (h_\varepsilon \bu_\varepsilon)+\partial_x(h_\varepsilon \bu_\varepsilon^2) \Bigg)  +\frac{1}{Fr^2}h_\varepsilon^{\beta} \partial_x(h_\varepsilon)=\varepsilon \Bigg(\frac{4}{{ R_e}}\partial_x(h_\varepsilon\partial_x \bu_\varepsilon) \Bigg) + \frac{1}{{ W_e}}h_\varepsilon\partial_x^3 h_\varepsilon 
- \alpha \frac{h_\varepsilon^2 \overline u_\varepsilon }{\he^n},
\end{cases}
\end{equation}
where $\beta+n \in (1,2)$.
The energy and BD entropy of system (\ref{swgg}) are given respectively by
\begin{equation}
\begin{aligned}
&\nonumber  \frac{1}{2} \left( \into \ep \he(x,T) \bu_\varepsilon^2(x,T)+\frac{1}{Fr^2}\frac{h_\varepsilon(x,T)^{\beta+1}}{\beta(\beta+1)}+\frac{(\partial_x h_\varepsilon(x,T))^2}{W_e}\right) ~dx +\int_{0}^{T}\into\frac{4\varepsilon }{R_e}h_\varepsilon (\partial_x \bu_\varepsilon)^2+  \alpha \frac{h_\varepsilon^2 \overline u_\varepsilon^2}{\he^n} ~dx~dt\\
& \nonumber = \frac{1}{2}\into \ep \frac{(m_0^{\ep})^2}{h_0^{\ep}} + \frac{1}{Fr^2}\frac{(h_0^{\ep})^{\beta+1}}{\beta(\beta+1)} + \frac{(\partial_x h_0^{\ep})^2}{W_e}~dx.
\end{aligned}
\end{equation}
and 
\begin{equation}
\begin{aligned}
&\nonumber  \frac{\ep}{2} \into \he(x,T)\ve(x,T)^2 - \he(x,T)\ue(x,t)^2~dx + \frac{4}{Re}\big[\int_{0}^{T}\into \frac{1}{Fr^2}\he^{\beta-1}( \partial_x\he)^2 ~dx+ \frac{1}{W_e} (\partial_{xx} \he)^2~dx~dt \big]\\
\nonumber 
&+ \frac{4}{Re} \into G_0(\he(x,T))~dx - \frac{4\ep}{Re}\int_{0}^{T} \into \he(\partial_x\ue)^2~dx~dt
\\
\nonumber 
&= \frac{\ep}{2} \into h_0^{\ep}\ve(x,0)^2 -\frac{m_0^{\ep}}{h_0^{\ep}}~dx 
+ \frac{4}{Re}  \into G_0(h_0^{\ep}(x))~dx.
\end{aligned}
\end{equation}
  Under the assumption of the convergence results, and choosing $m=\beta+n \in (1,2)$, we get that the BD entropy degenerates as well to  the BF dissipative entropy of system (\ref{lubg}) given by (\ref{bfg}). Remark that the link between the BD entropy and the BF entropy may also
be done in higher dimensions, see \cite{Gr} for BF entropy and \cite{BrDe2} for BD entropy. This
could help to perform the analysis in the bi-dimensional setting. 

\section{Mathematical results obtained using the link between BD and BF entropies}
    In this part, we aim at proving the existence of a  global in time weak solution for the lubrication model by passing to the limit in the viscous shallow water model with two different choices of the drag term, corresponding to two weights: $F(h)= h^2$ which results in a linear drag term, and $F(h)= h^2+ h^3$ which yields a nonlinear drag term. The latter weight has been used by A.L.Bertozzi in the physical and mathematical justification of the lubrication model \cite{Bert1}.
    The main theorem states:
\begin{theorem}\label{thm1}
	Given a sequence ${(\he,\ue)}_{\varepsilon}$ a global weak solution of (\ref{swg}), where $ h_0^{\ep} \ge 0$, then there exists a subsequence of $(\he,\ue)$  such that $(\he,\ue)$ converges to (h,u) a global weak solution of the lubrication system (\ref{lub1}) satisfying $h\ge0$, and the initial condition $ h|_{t=0}=h_0$, where $h_0$ is the weak limit of $h^\varepsilon_0$ in $H^1(\Omega$).
\end{theorem}
\vspace*{3mm}
The proof of the limit of the viscous shallow water model into a lubrication model is summarized in the following steps:
\begin{enumerate}
	\item  Assuming that (\ref{swg}) possesses a weak solution $(\he, \ue)$, bring first the physical energy and BD-entropy estimates  to get uniform bounds of the system's unknowns and thus get weak convergence up to a subsequence of these terms. Such solution has been constructed by D. Bresch and B. Desjardins in~\cite{BD}. See also \cite{LiLiXi} for more general interesting studies related to $1D$ compressible Navier-Stokes. 
	\item Use compactness theory to obtain strong convergence (mainly for $\he$).
	\item Pass the limit in the weak formulation of (\ref{swg}) to obtain that the solution is a weak solution of the lubrication theory.
\end{enumerate}
 It is important to remark that some studies have already analyzed the limit process from shallow-water to lubrication systems but with special pressure terms, see the nice papers \cite{MuWaWi}, \cite{KiLaNi} and recently \cite{FoKiTa}. See also the recent paper \cite{BrDu} where dissipative systems may be obtained from shallow-water type system through a quadratic change of time and no need of {\it a priori} drag terms.
 Note that in \cite{KiLaNi} and \cite{FoKiTa}, BD and BF entropies are interconnected for $F(h)=h$. 
   In \cite{ImMe}, we can find a result concerning the global existence of non-negative solutions for electrified thin films. Such systems contain nonlocal terms. As an example consider the following system
$$\partial_t h + \partial_x (h^3\partial_x(\partial_x^2 h  -  I(h)))
     = 0 \hbox{ in } \Omega = (0,1)$$
where $I(h)$ is a non-local elliptic operator of order 1 given by
 $$I(h) = \int_\Omega (h(y)-h(x)) \nu(x,y) \, dy$$
where for all $x,y\in \Omega$
$$\nu(x,y) = \frac{\pi}{2} \bigl(\frac{1}{1- \cos(\pi(x-y))} + \frac{1}{1- \cos(\pi(x+y))}   \bigr).$$
The system is supplemented by the following boundary and initial conditions
$$\partial_x h = h^3 \partial_x (\partial_x^2 h - I(h))= 0 \hbox{ on } \partial \Omega, \qquad
     h\vert_{t=0} = h_0 \hbox{ for } x\in \Omega.
$$ 
It is interesting to see that in order to construct  a solution to this lubrication equation, 
one can consider  the following shallow-water model
\begin{equation}\label{swnl}
\begin{cases}
\partial_t h_\varepsilon+\partial_x (h_\varepsilon\bu_\varepsilon)=0,\\
\varepsilon \Bigg(\partial_t (h_\varepsilon \bu_\varepsilon)+\partial_x(h_\varepsilon \bu_\varepsilon^2) \Bigg)  + h_\varepsilon \partial_x I(h_\varepsilon)=\varepsilon \partial_x(h_\varepsilon\partial_x \bu_\varepsilon) + h_\varepsilon\partial_x^3 h_\varepsilon 
- \alpha \frac{ \overline u_\varepsilon }{\he},
\end{cases}
\end{equation}
with appropriate boundary conditions.
Then using the energy estimate and the BD entropy, we can pass to the limit and
get the global existence of non-negative solutions of the system studied in
 \cite{ImMe}.  Note that compressible Navier-Stokes system with constant viscosities 
 and non-local term has been studied recently in \cite{CaWrZa} together  with the 
 long-time behavior of its solutions. The details will be given in the forthcoming paper \cite{BrCoMsNoSo}.





\end{document}